\documentclass[11pt]{article}
\usepackage[utf8]{inputenc}
\usepackage{amssymb}
\usepackage{amsmath}
\usepackage{amsthm}
\usepackage{color}
\usepackage{enumitem}
\usepackage{hyperref}
\usepackage[margin=1in]{geometry}
\usepackage{tikz}
\usepackage{bbm}

\theoremstyle{definition}
\newtheorem{definition}{Definition}[section]

\newtheorem{remark}[definition]{Remark}
\theoremstyle{plain}
\newtheorem{thm}[definition]{Theorem}
\newtheorem{prop}[definition]{Proposition}
\newtheorem{lemma}[definition]{Lemma}

\newtheorem{coro}[definition]{Corollary}

\newcommand{\vol}{\text{Vol}}
\newcommand{\sys}{\text{sys}}
\newcommand{\area}{\text{area}}

\title{Sharp systolic inequalities for invariant tight contact forms on principal $\mathbb{S}^1$-bundles over $\mathbb{S}^2$}
\author{Simon Vialaret\thanks{Université Paris-Saclay, CNRS, Laboratoire de mathématiques d'Orsay, 91405, Orsay, France \\simon.vialaret@universite-paris-saclay.fr}}
\date{}

\begin{document}

\maketitle

\begin{abstract}
    The systole of a contact form $\alpha$ is defined as the shortest period of closed Reeb orbits of $\alpha$. Given a non-trivial $\mathbb S^1$-principal bundle over $\mathbb S^2$ with total space $M$, we prove a sharp systolic inequality for the class of tight contact form on $M$ invariant under the $\mathbb S^1$-action. This inequality exhibits a behavior which depends on the Euler class of the bundle in a subtle way. As applications, we prove a sharp systolic inequality for rotationally symmetric Finsler metrics on $\mathbb S^2$, a systolic inequality for the shortest contractible closed Reeb orbit, and a particular case of a conjecture by Viterbo.
\end{abstract}

\section{Introduction}

A one-form $\alpha$ on a three-manifold $M$ is called a contact form when $\alpha \wedge \mathrm{d} \alpha$ is nowhere vanishing. It induces on $M$ the contact structure $\ker \alpha$, as well as a canonical vector field, namely the Reeb vector field $R$, defined by the following equations:
\begin{equation*}
	\begin{cases}
		\iota_R \alpha = 1\\
		\iota_R \mathrm{d}\alpha = 0.
	\end{cases}
\end{equation*}
The flow $(\phi^t)_{t\in \mathbb{R}}$ of $R$ is called the Reeb flow of $\alpha$. An important motivation for the study of Reeb flows comes from the fact that they generalize geodesic flows. Indeed, the cotangent bundle of any manifold can be equipped with a canonical one-form, called the Liouville form, whose restriction to the unit cosphere bundle induced by a Riemannian (or Finsler) metric on the base manifold is a contact form. The associated Reeb flow coincides with the cogeodesic flow of the metric. Moreover, the contact volume of this contact manifold is, up to a dimensional factor, the metric volume of the base manifold.

\begin{remark}
    If the context does not specify which contact form is being used, it will appear as an exponent on objects that depend on it.
\end{remark}

\subsection{The systolic ratio of a contact form} % (fold)
\label{par:systolic_ratio_of_a_contact_form}

Given a compact Riemannian manifold $(M,g)$ of dimension $n$, we define the systole of $g$, which we write $\sys(g)$, as the length of the shortest closed geodesic of $g$. On non-simply connected manifolds, one usually considers the homotopical systole, namely the length of the shortest non-contractible closed curve (necessarily a closed geodesic), but here we allow also contractible closed geodesics. The systolic ratio of $g$ is defined as $\rho(g) = \frac{\sys(g)^n}{\vol(g)}$, where $\vol(g)$ is the Riemannian volume of $M$. The study of the systole and of the systolic ratio is the objective of systolic geometry (see for example \cite{Katz_2007} and \cite{benedetti2021steps} for an introduction to the subject). Among objects of interest in systolic geometry is the supremum of $\rho(g)$ when $g$ varies in the set of all Riemannian metrics on a given manifold. Computing the possible values of this supremum, understanding how they depend on the topology of the manifold and describing the metrics for which the supremum is attained are important and difficult questions.

Motivated by the Riemannian setting, the authors of \cite{paiva_2014_contact} defined the systole of a contact form $\alpha$ on a manifold $M$, which we write $\sys(\alpha)$, as the minimal period of a periodic orbit of the Reeb flow of $\alpha$. If $M$ is a closed orientable manifold, this number is finite since Taubes proved in \cite{10.2140/gt.2007.11.2117} that Reeb flows on compact three-manifolds have a periodic orbit. The sytolic ratio of a contact form $\alpha$ on a contact manifold of dimension $2n+1$ is defined as the ratio $\rho(\alpha) = \frac{(\sys \alpha)^{n+1}}{\vol(\alpha)}$, where $\vol(\alpha)$ is the contact volume of $M$. Note that if $\alpha$ is the contact form induced on a unit cotangent bundle by a Riemannian metric, the metric and contact notions of systolic ratio coincide up to a dimensional factor.

In general, the contact systolic ratio behaves quite differently than its Riemannian counterpart. For example, write $C(\xi)$ for the set of contact forms inducing a given co-orientable contact structure $\xi$ on $M$. It was proved that the contact systolic ratio is unbounded on $C(\xi)$, first when $M$ is the three-sphere and $\xi$ is the unique tight contact structure on $\mathbb{S}^3$ in \cite{Abbondandolo2018}, on compact three-manifolds in \cite{Abbondandolo2019}, and in higher dimensions in \cite{Saglam_2021}. This has to be put in contrast with the fact that the Riemannian systolic ratio is known to be bounded on a large class of manifolds. Indeed, Gromov proved in \cite{MR0697984} that the homotopical systolic ratio is bounded from above by a constant depending only on the dimension whenever the non-simply connected manifold is essential. The boundedness of the Riemannian systolic ratio for $\mathbb{S}^2$ was proved by Croke in \cite{MR0918453}. There are no examples of compact manifolds for which the Riemannian systolic ratio is known to be unbounded. The study of the local maximizers of the systolic ratio on $C(\xi)$ was started in \cite{paiva_2014_contact}, and completed in \cite{Abbondandolo2018} (on $\mathbb{S}^3$ for the standard contact structure), in \cite{benedetti_kang_2020} (in dimension 3) and in \cite{Abbondandolo_Benedetti_2023} (in any dimension), where it was proved that the local maxima of the systolic ratio are exactly Zoll contact forms, namely contact forms whose Reeb flow only has periodic orbits, all of them with same minimal period. Again, this has to be compared with the metric situation where maximizers of the systolic ratio are poorly understood.
% subsection systolic_ratio_of_a_contact_form (end)

\subsection{Invariant tight contact forms and the main result} % (fold)

\paragraph{Notations} % (fold)
\label{par:notations}

% paragraph notations (end)
Let $M$ be a compact oriented three-manifold, endowed with an $\mathbb S^1$-principal bundle structure $M \overset{\pi}{\rightarrow} \mathbb S^2$, where $\mathbb{S}^1$ denotes $\mathbb{R}/\mathbb{Z}$. Every such manifold is topologically either $\mathbb S^1 \times \mathbb S^2$ when the bundle is trivial, or a lens space $L(n,1)$, which is a quotient of $\mathbb{S}^3$ seen as the unit sphere in $\mathbb C^2$ by the action of $\mathbb Z / n \mathbb Z$ generated by $(z_1, z_2) \mapsto (e^{\frac{2 \pi i}{n}} z_1, e^{\frac{2 \pi i}{n}} z_2)$. For $x \in M$, the orbit of $x$ under the $\mathbb S^1$-action will be called the $\mathbb{S}^1$-orbit of $x$, and denoted by $\mathbb{S}^1\cdot x$. We will write $\Theta$ for the vector field generating the action, meaning that $\Theta_x = \left.\frac{d}{d t} \right\vert_{t = 0} t \cdot x$. A contact form on $M$ is said to be positive if the volume form $\alpha \wedge \mathrm{d}\alpha$ induces the orientation of $M$. Finally, we say that $\alpha$ is an $\mathbb{S}^1$-invariant contact form on $M$ if $\alpha$ is a contact form satisfying $\mathcal{L}_\Theta \alpha = 0$.

$\mathbb{S}^1$-principal bundles over a base $B$ are classified by homotopy classes of maps from $B$ to the classifying space of $\mathbb{S}^1$, which is $\mathbb{C}P^\infty$. This allows us to define the Euler class of an $\mathbb{S}^1$-principal bundle as the pullback of $-1 \in H^2(\mathbb{C}P^\infty, \mathbb{Z})$ (recall that $\mathbb{C}P^\infty$ is naturally oriented) by the associated map. Since $M$ is oriented, the principal bundle structure induces an orientation on the base $\mathbb{S}^2$ as follows: if $(X,Y; \Theta)$ is a positive local frame on $M$, then $(\pi_* X,\pi_* Y)$ is a local positive frame on $\mathbb S^2$. We then define the Euler number $e$ of the bundle as the integer given by pairing the fundamental class of $\mathbb{S}^2$ with the Euler class. The data of an orientation on $\mathbb{S}^2$ and of an integer $e$ determines a cohomology class in $H^2(\mathbb{S}^2, \mathbb{Z})$, which by the classification of $\mathbb{S}^1$-principal bundles determines a unique $\mathbb{S}^1$-principal bundle over $\mathbb{S}^2$ with Euler number $e$ and canonically oriented total space. We give a geometric description of the Euler number in Remark \ref{remark:euler_number}.

In view of the fundamental dichotomy between tight and overtwisted contact structures in contact topology (see \cite{geiges_2008}), we focus on contact forms inducing a tight contact structure in this paper. We write $\Omega(e)$ for the set of $\mathbb{S}^1$-invariant positive tight contact forms on the total space $M$ of the above $\mathbb{S}^1$-principal bundle with Euler number $e$. In \cite{Giroux1999StructuresDC}, Giroux proved the following results: 
\begin{itemize}
    \item If $e < 0$, then the contact forms in $\Omega(e)$ induce the unique (up to isotopy) positive contact structure that is transverse to the orbits of the action;
    \item if $e \geq 0$, then the contact forms in $\Omega(e)$ induce a contact structure whose set of Legendrian $\mathbb S^1$-orbits is non-empty and connected.
\end{itemize}

The classification of invariant contact structures on principal $\mathbb S^1$-bundles have been achieved in \cite{Lutz1977} and \cite{Giroux1999StructuresDC}, see also \cite[Section 5]{GeigesHedickeSaglam} for a detailed discussion of $\mathbb S^1$-invariant contact structures on $\mathbb S^3$. We give below examples of tight $\mathbb{S}^1$-invariant contact forms on bundles over $\mathbb S^2$.

\begin{itemize}
    \item The standard contact form $\alpha_0 = \frac{1}{4i}\sum\limits_{i=1}^{2}\left(\bar{z_i} \mathrm{d}z_i - z_i \mathrm{d}\bar{z_i} \right)$ on $\mathbb{S}^3 \subset \mathbb{C}^2$ is invariant under the $\mathbb{S}^1$-action $t \cdot (z_1, z_2) \mapsto (e^{2 \pi i t}z_1, e^{2 \pi i t}z_2)$ given by its Reeb flow. The corresponding $\mathbb{S}^1$-bundle has Euler number $-1$. More generally, any positive Zoll contact form on a manifold $M$ is invariant under the $\mathbb{S}^1$-action induced by its Reeb flow, and the induced principal bundle structure on $M$ has negative Euler number. Zoll contact forms on non-trivial $\mathbb S^1$-principal bundles over a symplectic manifold with integral symplectic form can be obtained by the Boothby--Wang construction, see Section 7.2 in \cite{geiges_2008}.
    
    \item Another class of examples come from a class of metrics on $\mathbb{S}^2$. A Finsler metric on a manifold $M$ is a continuous non-negative function $F \colon T M \rightarrow \mathbb{R}^+$ satisfying the following properties:
    \begin{itemize}
        \item $F$ is smooth on the complement of the zero-section;
        \item if $u \in TM$ is not zero then $\mathrm{d}^2F(u)(w,w) > 0$ for all $w \in T^{\text{vert}}_u TM \setminus (\mathbb R u)$;
        \item $F(\lambda v) = \lambda F(v)$ for $\lambda \geq 0$;
        \item $F(v) > 0$ if $v \neq 0$.
    \end{itemize}
    Riemannian metrics are special cases of Finsler metrics, where $F$ is the Euclidean norm induced by the metric. As in the Riemannian setting, a Finsler metric induces a contact form $\alpha_F$ on its unit tangent bundle. We call a Finsler metric on $\mathbb{S}^2$ a metric of revolution if $F$ is invariant under a circle action on $\mathbb{S}^2$ that has exactly two fixed points and is free on their complement. This action lifts to a free $\mathbb{S}^1$-action on $T^1 \mathbb{S}^2$ the unit tangent bundle for $F$. The induced principle bundle structure induced on $T^1 \mathbb{S}^2$ has Euler number 2. The contact form $\alpha_F$ is invariant under that action. Hence, $\alpha_F$ is in $\Omega(2)$.

    \item If $\alpha$ is a positive $\mathbb{S}^1$-invariant contact form on $M$ and $f \colon M \rightarrow ]0, +\infty[$ is invariant under the same action, then $f \alpha$ is an $\mathbb{S}^1$-invariant positive contact form as well.

\end{itemize}

The main result of this work is a sharp systolic inequality for invariant tight contact forms on non-trivial $\mathbb S^1$-principal bundles over $\mathbb S^2$. This inequality exhibits a behavior which depends on the Euler number of the bundle in a non-trivial way.

\begin{thm}
    \label{thm:ineg_sys_tight}
    \begin{itemize}
        \item If $e < 0$ or if $e \in \{1,2\}$, then for all $\alpha$ in $\Omega(e)$, we have
        	\[
        		\sys(\alpha)^2 \leq \frac{1}{|e|} \vol(\alpha).
        	\]
        	Moreover, equality holds if and only if $\alpha$ is Zoll.
        \item If $e > 2$, then for all $\alpha$ in $\Omega(e)$, we have
        	\[
        		\sys(\alpha)^2 < \frac{1}{2} \vol(\alpha).
        	\]
            Moreover, the supremum of the systolic ratio $\rho(\alpha)$ on $\Omega(e)$ is $\frac{1}{2}$ (and is not achieved).
    \end{itemize}
	
\end{thm}

\begin{remark}
    \begin{itemize}
        \item Since the set of invariant contact forms with fixed volume has infinite diameter for the Banach--Mazur pseudo-metrics on the set of contact forms introduced in \cite{RosenZhang} and \cite{abbondandolo2023symplectic}, this result shows that boundedness of the systolic ratio is not a purely metric phenomenon, as the dichotomy ``usually bounded metric systolic ratio''/``never bounded contact systolic ratio'' might suggest. 
        
        \item We focus here on tight contact forms, for which we can prove the sharp systolic inequality above. However, the existence of a systolic inequality for invariant contact forms is not a purely tight phenomenon. We prove in a forthcoming paper a general systolic inequality for invariant contact forms on $\mathbb S^1$-principal bundles in dimension 3, valid for both tight and overtwisted invariant contact forms.

        \item The fact that the supremum of $\rho(\alpha)$ on $\Omega(e)$ is not achieved when $e > 2$ can be understood as follows\footnote{The author wishes to thank Christian Lange for pointing out that the systolic inequality had to be different for Euler numbers greater than 2.}: the classification of Zoll contact forms in dimension 3, achieved in \cite{benedetti_kang_2020}, implies that there is no positive Zoll contact form on $M$ when $e>2$. Moreover, contact forms that are local maximizers of the systolic ratio $\rho(\alpha)$ in the class of $\mathbb{S}^1$-invariant contact forms are necessarily Zoll, as proved in \cite[Theorem 3.4]{paiva_2014_contact}. Thus the equality cannot be reached for the systolic inequality for invariant forms on $M$ when $e > 2$.

        \item The statement of Theorem \ref{thm:ineg_sys_tight} leaves the case of invariant contact forms on $\mathbb S^1 \times \mathbb{S}^2$ open. This will be an objective of further work.
    \end{itemize}
\end{remark}

% subsection a_systolic_inequality_for_invariant_contact_forms (end)

\subsection{Applications}
\subsubsection{Systolic inequality for Finsler metrics of revolution on \texorpdfstring{$\mathbb S^2$}{S2}}

\label{ssub:Finsler}
As discussed in \cite{MR0918453}, the supremum of the Riemannian systolic ratio for $\mathbb{S}^2$ is finite, larger than the one of the round sphere but its value is not known. The authors of \cite{ABHS_revol} proved that the systolic ratio of a Riemannian sphere of revolution is maximized by Zoll metrics. They also extended this result to a special class of $\mathbb{S}^1$-symmetric Finsler metrics. We now explain how Theorem \ref{thm:ineg_sys_tight} allows us to extend this result to all $\mathbb{S}^1$-invariant Finsler metrics on $\mathbb{S}^2$.

A subtlety arising in Finsler geometry is that, as opposed to the Riemannian setting, there is no canonical notion of volume. As a matter of fact, several notions of volume coexist in the literature. As it is of symplectic nature, we will consider the Holmes--Thomson volume with the normalization convention of \cite{ABHS_revol}: if $(S,F)$ is a Finsler surface, and $U$ an open subset of $S$, we define the Holmes--Thomson area of $U$ as
$$
    \area_{HT}(U) = \frac{1}{2 \pi} \vol(T^{*1}U, \lambda_{\text{Liouville}}).
$$

In other words, $\area_{HT}(U)$ is, up to a normalizing factor, the contact volume of the unit cotangent bundle over $U$ for the canonical Liouville form. We will write $\area_{HT}(F)$ for the Holmes--Thomson area of a Finsler two-sphere $(\mathbb S^2, F)$. More details on volumes in Finsler geometry can be found in \cite{paiva2004volumes}. We show that the following systolic inequality holds for any Finsler metric of revolution.
\begin{coro}
    \label{coro:finsler}
    If $F$ is a Finsler metric of revolution on $\mathbb S^2$, then the inequality
    $$
        \sys(F)^2 \leq \pi \area_{HT}(F)
    $$
    holds, with equality if and only if $F$ is a Zoll metric.
\end{coro}

\begin{proof}
    Let $F$ be a Finsler metric of revolution on $\mathbb S^2$. The contact form $\alpha_F$ is in $\Omega(2)$. Theorem \ref{thm:ineg_sys_tight} applied to $\alpha_F$ gives
    \[
        \sys(\alpha_F)^2 \leq \frac{1}{2} \vol(\alpha_F).
    \]
    Since $\sys(F) = \sys(\alpha_F)$ and $\vol(\alpha_F) = 2 \pi \area_{HT}(F)$, we get the desired inequality, with equality if and only if $F$ is Zoll.
\end{proof}

\subsubsection{Contractible systole and Besse contact forms} % (fold)
\label{ssub:contractible_systole_and_besse_contact_forms}

Another class of special contact forms, generalizing the Zoll ones, is the class of Besse contact forms. A contact form is Besse when its Reeb flow only has periodic orbits (not necessarily of same minimal period as in the Zoll case). A classical result due to Wadsley \cite{Wadsley_1975}, or in the three dimensional case an earlier result of Epstein \cite{Epstein}, implies that the Reeb flow of a Besse contact form is periodic. The theorem of Epstein implies that all orbits but finitely many ones of the Reeb flow of a Besse form have the same minimal period $T$, and the minimal period of the finitely many others divides $T$. Similarly to the principle conjectured in \cite{paiva_2014_contact} and proved in \cite{Abbondandolo2018}, \cite{benedetti_kang_2020}, \cite{Abbondandolo_Benedetti_2023}, that Zoll contact forms are exactly the local maximizers of the systolic ratio, some variations of the systolic ratio were studied in \cite{abbondandolo:ensl-03357629}, \cite{baracco2023local} and \cite{MR4529024} for which Besse contact forms are exactly the local maximizers. In \cite{MR4529024}, the authors study the contractible systole of a contact form $\alpha$, denoted $\sys_{\text{contr}}(\alpha)$, namely the shortest period of a contractible periodic orbit of the Reeb flow of $\alpha$. They prove a systolic inequality for the contractible systole for the contact forms induced by a class of singular Riemannian metrics on $\mathbb S^2$ called rotationally symmetric spindle orbifolds (Theorem A in \cite{MR4529024}), generalizing the main result of \cite{ABHS_revol}. Moreover, equality holds exactly when the spindle orbifold is Besse.

Theorem \ref{thm:ineg_sys_tight} also implies a systolic inequality for the contractible systole, stated in Corollary \ref{coro:sys_contr} below, whose proof will be given in \ref{sub:proof_coro_contractible}. It recovers as a particular case Theorem A in \cite{MR4529024}. Depending on the Euler number of the circle action, the maximizers of the contractible systolic ratio are either Zoll contact forms, or Besse contact forms with exactly two singular orbits, which have the same minimal period. When $e>2$, Besse contact forms with two singular orbits of same minimal period can be obtained as a quotient of the standard contact form $\alpha_0$ on $\mathbb{S}^3 \subset \mathbb{C}^2$ via the $\mathbb{Z}/ e \mathbb{Z}$ action generated by $(z_1, z_2) \mapsto (e^{\frac{2 \pi i}{e}} z_1, e^{-\frac{2 \pi i}{e}} z_2)$. If $T$ is the minimal period of the Reeb flow $\alpha_0$, the contact form induced on the quotient has two singular orbits of period $\frac{T}{e}$. The other Reeb orbits have period $T$ if $e$ is odd, $\frac{T}{2}$ if $e$ is even. Alternatively, those Besse contact forms can be described as spindle orbifolds as in \cite{MR4529024}.

\begin{coro}
    \label{coro:sys_contr}
    Let $e \neq 0$, then for all $\alpha \in \Omega(e)$,
    \[
        \sys_{\text{contr}}(\alpha)^2 \leq |e| \vol(\alpha).
    \]
    Moreover, equality holds if and only if $\alpha$ is Zoll when $e \leq 2$, and if and only if $\alpha$ is Besse with two singular orbits of same minimal period, when $e>2$. In the latter case, the singular orbits have period $\frac{T}{e}$ if $e$ is odd, $\frac{2T}{e}$ is $e$ is even, where $T$ is the minimal period of the regular orbits.
\end{coro}

% subsubsection contractible_systolic_ratio_and_besse_contact_forms (end)

\subsubsection{Viterbo conjecture for invariant domains}
\label{ssub:viterbo}
Another corollary of Theorem \ref{thm:ineg_sys_tight} is an application to a conjecture of Viterbo. Let $\omega_0 = \sum\limits_{i=1}^n \mathrm{d}x_i \wedge \mathrm{d}y_i$ be the standard symplectic form on $\mathbb R^{2n}$. We write $\vol(U, \omega_0^n)$ for the symplectic volume of a set $U \subset \mathbb{R}^{2n}$. In \cite{Viterbo_2000} (see also \cite{MR4413022} for an overview of variations of this conjecture), Viterbo conjectures that for any symplectic capacity $c$ and any open bounded convex domain $U \subset \mathbb{R}^{2n}$ the inequality
\begin{equation}
    \label{eq:conj_viterbo}
    c(U)^n \leq \vol(U, \omega_0^n)
\end{equation}
holds, and that equality holds if and only if $U$ is symplectomorphic to a ball. He also proves the inequality up to a dimensional factor. Since any convex domain $U$ in $\mathbb{R}^{2n}$ is of contact type, its boundary carries a contact form $\alpha_U$. For many normalized capacities $c$ it is known that $c(U)=\sys(\alpha_U)$ whenever $U$ is a smooth bounded convex domain. This is for example the case for the capacities $c_{HZ}, c_{SH}, c_1^{EH}, c_1^{CH}, c_1^{Alt}$, following the work of Ekeland, Hofer, Zehnder, Abbondandolo, Kang, Irie and Hutchings (\cite{hryniewicz2023hopf} provides a related discussion). By a consolidated tradition, we denote their common value on such domains by $c_{EHZ}$. As a consequence, systolic inequalities for convex domains can be seen as special cases of Viterbo conjecture.

It is now known that Viterbo conjecture does not hold at this level of generality. Indeed, in \cite{haimkislev2024counterexample}, the authors constructed explicit examples of bounded convex domains $U \subset \mathbb{R}^4$ satisfying $\frac{c_{EHZ}(U)^2}{\vol(U, \omega_0^2)} = \frac{\sqrt{5}+3}{5} > 1$, as well as examples in any higher dimension, disproving Viterbo conjecture. On the other hand, Viterbo conjecture has been confirmed for several subclasses of convex domains (see \cite{haimkislev2024counterexample} and the references therein). Furthermore, in view of its connection with Mahler's conjecture in convex geometry (see \cite{10.1215/00127094-2794999}), it is interesting to explore the range of its validity. Theorem \ref{thm:ineg_sys_tight} allows us to show that Viterbo conjecture holds in explicit situations, some of which are new.

Let $(\phi^t)_{t \in \mathbb{S}^1}$ be a symplectic circle action on $(\mathbb{R}^4, \omega_0)$, and $U$ a smooth compact domain that is invariant under this action and with contact-type boundary. The action restricts to the boundary of $\partial U$, which is diffeomorphic to $\mathbb{S}^3$. If $Y$ is a Liouville vector field on an invariant neighborhood of $\partial U$, transverse to $\partial U$, then averaging $Y$ along the orbits of the action produces an invariant Liouville vector field $\tilde Y$, which is transverse to $\partial U$. The contraction $\iota_{\tilde Y}\omega_0$ then induces a tight contact form $\alpha_U$ on $\partial U$. Furthermore, $\alpha_U$ is invariant under the circle action restricted to $\partial U$. Whenever this action happens to be free on $\partial U$, Theorem \ref{thm:ineg_sys_tight} implies that 
\[
    \sys(\alpha_U)^2 \leq \vol(\alpha_U) = \vol(U, \omega_0^n),
\]
with equality if and only if $\alpha_U$ is Zoll, which by Proposition 4.3 in \cite{Abbondandolo2018} is equivalent to the fact that $U$ is symplectomorphic to a ball. This leads us to the corollary below.

\begin{coro}
    \label{coro:viterbo_invariant}
    Let $(\phi^t)_{t \in \mathbb{S}^1}$ be a symplectic circle action on $(\mathbb{R}^4, \omega_0)$, and $U$ a smooth convex domain that is invariant under this action. If the restricted action on $\partial U$ is free, then we have
    \[
        c_{EHZ}(U)^2 \leq \vol(U, \omega_0^n),
    \]
    with equality if and only if $U$ is symplectomorphic to a ball.
\end{coro}

As a concrete example, let $\alpha_0 = \frac{1}{2}\sum \limits_{i=1}^2 x_i \mathrm{d}y_i - y_i \mathrm{d}x_i $ be the radial Liouville form on $\mathbb{R}^4$. Let $(\Phi_t^+)_{t \in \mathbb S^1}$ and $(\Phi_t^-)_{t \in \mathbb S^1}$ be the two circle actions on $\mathbb C^2$ defined by 
$$
    \Phi^+_t(z_1, z_2) = (e^{2i\pi t} z_1, e^{2i\pi t} z_2) \hspace{1cm} \text{ and} \hspace{1cm} \Phi^-_t(z_1, z_2) = (e^{2i\pi t} z_1, e^{-2i\pi t} z_2).
$$

Identifying each $(x_i,y_i)$-plane of $\mathbb R^4$ with a copy of $\mathbb C$ yields two symplectic actions on $\mathbb R^4$. Those two actions are not symplectomorphic, since they define two principal bundle structures on $\mathbb{S}^3$ with different Euler number, hence not conjugated via an orientation preserving diffeomorphism. Moreover, those actions are free on $\mathbb{R}^4 \setminus \{0\}$. Hence, as implied by Corollary \ref{coro:viterbo_invariant}, any convex domain $U$ that is invariant under any of those action satisfies the Viterbo conjecture. We have respectively $e=-1$ and $e=1$ for the induced actions of $\Phi^+$ and $\Phi^-$ on $\partial U$. As we will see from Proposition \ref{prop:ineg_sys_facile}, this makes Corollary \ref{coro:viterbo_invariant} essentially trivial for $\Phi^+$, but non-trivial and new for $\Phi^-$. The non-triviality of the result for $\Phi^-$ can be illustrated by the fact that the restriction of $\Phi^-$ to the unit sphere in $\mathbb C^2$ is a lift to $\mathbb S^3$ of the circle action on $T_1 \mathbb S^2 = \mathbb RP^3$ described in \ref{ssub:Finsler}, which follows from the classification of $\mathbb S^1$-principal bundles over $\mathbb S^2$ by the Euler number. As such, Corollary \ref{coro:finsler} can be seen as a consequence of Corollary \ref{coro:viterbo_invariant} for the action $\Phi^-$.

\subsection{A sketch of the proof}
The proof of Theorem \ref{thm:ineg_sys_tight} splits into two cases, depending on the sign of the Euler number of the bundle. When it is negative, the result is proved by an easy argument, see for example \cite{hutchings_2013}. This is the content of Proposition \ref{prop:ineg_sys_facile}.

When the bundle has positive Euler number and the invariant contact structure has a non-empty connected set of Legendrian $\mathbb S^1$-orbits, the proof goes as follows.
\begin{itemize}
    \item The first technical tool coming into play is a well-suited choice of surface of section for the Reeb flow of an invariant contact form. Given a flow on a manifold $M$, a surface of section for this flow is an embedded surface $S \subset M$ transverse to the flow, whose boundary is a union of orbits of the flow, and with the property that for any $x$ in $S$, the orbit of $x$ by the flow hits $S$ again in forward time.
    Knowing that such a surface exists allows us to reduce the study of the flow, or at least part of the flow, to the study of the first-return map defined on $S$, and hence to reduce the dimension of the problem by one. The use of surfaces of section in dynamics dates back to Poincaré in his study of the three-body problem, and has been very fruitful since then (see for example the surveys \cite{Ghys_2009} and \cite{SALOMAO_HRYNIEWICZ_2019}). It has proved to be helpful in the quantitative study of Reeb flows (see for example \cite{Albers2023}, \cite{Moreno_van_Koert_2022}, \cite{Chaidez2022}, \cite{Albach2021}, \cite{Abbondandolo2023}, \cite{Contreras2022}, \cite{colin2022generic}, \cite{contreras2022surfaces}, \cite{edtmair2022disklike} for recent applications) and in particular for proving sharp systolic inequalities in both the metric and contact setting (see for example \cite{Abbondandolo2017}, \cite{Abbondandolo2018}, \cite{Saglam_2021}, \cite{ABHS_revol}, \cite{benedetti_kang_2020}, \cite{Asselle2021}). We construct in Section \ref{sub:invariant_surface_of_section} an invariant surface of section $\Sigma$ for the Reeb flow of any $\mathbb S^1$-invariant contact form. This surface of section is a union of embedded cylinders, each of them invariant under the $\mathbb S^1$-action. Even though $\Sigma$ has a priori several connected components, we will only make use of the one that intersects the set of Legendrian $\mathbb S^1$-orbits.

    \item By $\mathbb{S}^1$-invariance of the surface of section and of the Reeb flow, studying the first-return map boils down to the study of a real-valued function, which we obtain using the decomposition of invariant contact forms obtained in Section \ref{sub:a_normal_form_for_invariant_contact_forms}. This function both encodes information on periodic orbits of the Reeb flow and their periods, and on the contact volume of $M$ via its integral. Building on techniques used in \cite{ABHS_revol}, we use this function to find a suitable periodic point of the first-return map, corresponding to a periodic Reeb orbit with controlled period. This yields the result.
\end{itemize}

\subsection{Acknowledgments}
The author is very grateful to Alberto Abbondandolo and Rémi Leclercq for numerous stimulating conversations and advice, as well as to Christian Lange for his careful reading of the first version of this work. The author is partially funded by the Deutsche Forschungsgemeinschaft (DFG, German Research Foundation) – Project-ID 281071066 – TRR 191 and by the ANR grant 21-CE40-0002 (CoSy).

\section{Generalities on \texorpdfstring{$\mathbb{S}^1$}{S1}-invariant contact forms}
\label{sec:generalites}
We introduce in this section the first integral of the Reeb flow induced by the $\mathbb{S}^1$-symmetry, and give some of its properties. It is a general fact that to a symmetry of a dynamical system corresponds an invariant of the dynamics. In our case, this invariant is given by the moment map of the action. Recall that $\Omega(e)$ is the set on $\mathbb{S}^1$-invariant positive tight contact forms on the total space of the principal $\mathbb{S}^1$-bundle over $\mathbb{S}^2$ with Euler number $e$.

\begin{definition}
	Let $\alpha \in \Omega(e)$. The moment map of $\alpha$ is the function $\tilde K \colon M \rightarrow \mathbb R$ defined by
	\[
		\tilde K = i_{\Theta} \alpha.
	\]
\end{definition}

By $\mathbb{S}^1$-invariance, $\tilde K$ induces a function on $\mathbb{S}^2$, which we denote $K$.

\begin{prop}
\label{prop:K_invariant}
	$\tilde K$ is left invariant by the Reeb flow.
\end{prop}

\begin{proof}
	We must prove that $\mathcal L_R \tilde K = 0$, which follows from Cartan's formula:
	\[
		\mathcal L_R \tilde K = (i_R \mathrm{d} + \mathrm{d}i_R)i_\Theta \alpha = i_R \mathrm{d} i_\Theta \alpha = i_R \mathcal L_\Theta \alpha - i_R i_\Theta \mathrm{d}\alpha.
	\]
	Since $\alpha$ is $\mathbb{S}^1$-invariant, one has $\mathcal L_\Theta \alpha = 0$, hence the first term is zero. The second term vanishes as $R$ is in the kernel of $\mathrm{d}\alpha$.
\end{proof}

An important property of the moment map is that there is a bijection between the critical points of $K$ and the Reeb orbits of $\alpha$ that coincide with an $\mathbb{S}^1$-orbit as sets. Moreover, the period of those Reeb orbits is up to a sign the corresponding critical value of $K$.

\begin{prop}
	\label{prop:crit_K_S1_Reeb}
    The differential of $\tilde K$ is $\mathrm{d}\tilde K = - i_\Theta \mathrm{d}\alpha$. Moreover, the critical points of $K$ are in bijection with $\mathbb{S}^1$-orbits of critical points of $\tilde K$, and themselves are exactly the $\mathbb{S}^1$-orbits that coincide as sets with Reeb orbits. The associated critical value of $K$ is the minimal period of this Reeb orbit up to sign, depending on whether the orientations of the corresponding Reeb orbit and $\mathbb{S}^1$-orbit coincide or not.
\end{prop}

\begin{proof}
	We can apply Cartan's formula to compute $\mathrm{d}\tilde K$:
	\[
		\mathrm{d}\tilde K = \mathrm{d}i_{\Theta}\alpha = \mathcal L_{\Theta} \alpha - i_{\Theta}\mathrm{d}\alpha = - i_{\Theta}\mathrm{d}\alpha
	\]
	which proves the first claim. We obtain from $\mathrm{d}K \circ \mathrm{d}\pi = \mathrm{d}\tilde K$ that $\mathrm{d}K$ vanishes at $x \in \mathbb{S}^2$ if and only if $\Theta_y$ is in the kernel of $\mathrm{d}\alpha_y$, with $\pi(y) = x$. Since the kernel of $\mathrm{d}\alpha$ is spanned by $R$, $\mathrm{d}K$ vanishes at $x$ if and only if $\Theta_y$ and $R_y$ are parallel. By $\mathbb{S}^1$-invariance of $R$ and $\Theta$, this happens if and only if $\mathbb{S}^1\cdot y$ and the Reeb orbit starting at $y$ coincide as sets. In that case, writing $\gamma$ for the Reeb orbit starting at $y$, with minimal period $T$, and $\gamma_\Theta$ for $\gamma$ parameterized by the flow of $\Theta$, one has 
	\[
		T = \int_\gamma \alpha = \left| \int_{\gamma_\Theta} \alpha\right| = \left| \int_{0}^{1} \alpha_{t.y}(\Theta) \mathrm{d}t\right| = |\tilde K(\gamma)| = |K(x)|.
	\]
\end{proof}

A corollary of Proposition \ref{prop:crit_K_S1_Reeb} is that $0$ is always a regular value of $K$. Indeed, $K$ vanishes exactly when $\Theta$ lies in the contact structure, hence $\Theta$ and $R$ cannot be parallel when $K$ vanishes. We finish this section by proving the inequality stated in Theorem \ref{thm:ineg_sys_tight} when the Euler number is negative. For a tight invariant contact structure, this happens exactly when the action of $\mathbb{S}^1$ is transverse to the contact structure. In that specific case, this inequality was already proved in the blog post \cite{hutchings_2013} of Hutchings (see also Proposition 1.4 in \cite{MR4413022}).

\begin{prop}
    \label{prop:ineg_sys_facile}
	Let $e<0$ and $\alpha \in \Omega(e)$. Then the inequality
	\[
		\sys(\alpha)^2 \leq \frac{\vol(\alpha)}{|e|}
	\]
	holds. Moreover, equality holds if and only if $\alpha$ is Zoll.
\end{prop}

\begin{proof}
	Since $e<0$, we know that the $\mathbb{S}^1$-action is transverse to the contact structure. As a consequence, $\tilde K^\alpha$ does not vanish. This means that $\alpha_0 = \frac{\alpha}{\tilde K^\alpha}$ is an $\mathbb{S}^1$-invariant contact form on $M$, and $\tilde K^{\alpha_0} = 1$. Since $\tilde K^{\alpha_0}$ is constant, all the Reeb orbits of $\alpha_0$ are also $\mathbb{S}^1$-orbits, and have period $1$. In other words, $\alpha_0$ is a Zoll contact form on $M$, and as such has volume $|e|$ (see Proposition 3.3 in \cite{paiva_2014_contact}). The volume of $\alpha$ can be written
	\[
		\vol(\alpha) = \int_{M} \alpha \wedge \mathrm{d}\alpha = \int_{M} (\tilde K^\alpha)^2 \alpha_0 \wedge \mathrm{d}\alpha_0.
	\]
	Since $\tilde K^\alpha$ does not vanish, $(\tilde K^\alpha)^2$ has a positive minimum that we write as $K_0^2$. This critical point of $\tilde K^\alpha$ corresponds to a closed Reeb orbit of period $|K_0|$, hence the systole of $\alpha$ is not larger than $|K_0|$. Finally, one has
	\[
		\vol(\alpha) = \int_{M} (\tilde K^\alpha)^2 \alpha_0 \wedge \mathrm{d}\alpha_0 \geq K_0^2 \int_{M} \alpha_0 \wedge \mathrm{d}\alpha_0 = |e|K_0^2 \geq |e|\sys(\alpha)^2.
	\]
 Equality holds if and only $\tilde K^\alpha$ is constant, i.e. if and only if $\alpha$ is Zoll.
\end{proof}

\begin{remark}
    \label{remark:tight_K_connected}
    From now on, we will assume that $e>0$. The Giroux classification of tight invariant contact structure implies that for $\alpha$ in $\Omega(e)$, the set of Legendrian $\mathbb S^1$-orbits is non-empty and connected, in other words  the regular level set $(K^\alpha)^{-1}(0)$ is not empty and connected, and therefore a circle.
\end{remark}

\section{Invariant surfaces of section and potentials}
\label{sec:inv_surface_of_section_and_potentials}
This section contains the two main geometrical constructions underlying the proof of Theorem \ref{thm:ineg_sys_tight}, namely an invariant surface of section, and a decomposition for invariant contact forms. These two constructions will allow us to introduce a family of potentials encoding most of the dynamic of the Reeb flow. We end this section by giving some important properties of these potentials. 

\subsection{An invariant surface of section}
\label{sub:invariant_surface_of_section}
Let us fix $\alpha$ in $\Omega(e)$. In this section, we will construct an $\mathbb{S}^1$-invariant surface of section for the Reeb flow of $\alpha$. We write $Z$ for the regular level set $\tilde K^{-1}(0) \subset M$. Let $\mathcal C \subset \mathbb{S}^2$ be the union of the connected components of levels sets of $K$ that contain a critical point. It follows from the fact that $\mathrm{d}K$ does not vanish on $\mathbb{S}^2 \setminus \mathcal C$, that $\mathbb{S}^2 \setminus \mathcal C$ decomposes as a union of cylinders $(\Gamma_i)_{i \in I}$. By Remark \ref{remark:tight_K_connected}, $\pi(Z)$ is connected. Hence, there is a unique index $i_0$ such that $K$ has zeroes on $\Gamma_{i_0}$.

For each $i \in I$, we write $\epsilon_i^+$ and $\epsilon_i^-$ for the two components of $\partial \Gamma_i$, with $K(\epsilon_i^+) > K(\epsilon_i^-)$, and $K_i^+$ (resp. $K_i^-$) for $K(\epsilon_i^+)$ (resp. $K(\epsilon_i^-)$). We choose $x_i^+$ (resp. $x_i^-$) to be an arbitrary critical point of $K$ in $\epsilon_i^+$ (resp. $\epsilon_i^-$). Finally, we choose a curve $\gamma_i \colon ]0,1[ \rightarrow \Gamma_i$ such that:
\begin{itemize}
    \item $\gamma_i$ is transverse to the level sets of $K$;
    \item $\lim_{s \rightarrow 0}\gamma_i(s) = x_i^-$ and $\lim_{s \rightarrow 1}\gamma_i(s) = x_i^+$.
\end{itemize}

Writing $\tilde \gamma_i$ for a lift of $\gamma_i$ to $M$, we set $\Sigma_i$ as the image of $\tilde \gamma_i$ under the $\mathbb{S}^1$-action, i.e. $\Sigma_i = \mathbb{S}^1\cdot \tilde \gamma_i = \pi^{-1}(\gamma_i)$. The $\Sigma_i$ are disjoint cylinders, whose boundary components have the form $\mathbb{S}^1 \cdot x$, with $x$ a critical point of $K$, and as such are both closed Reeb orbits and $\mathbb{S}^1$-orbits. We write $\Sigma$ for their union $\bigcup \limits_{i \in I} \Sigma_i$. The following proposition states that $\Sigma$ is a surface of section.

\begin{prop}
	\label{prop:temps_de_retour_fini}
    The vector field $R$ is transverse to $\Sigma$. Moreover, for all $x$ in $\Sigma$, the Reeb orbit of $x$ intersects $\Sigma$ in positive time.
\end{prop}

\begin{proof}
    Let $x \in \Sigma$. Since the vector field $R$ is nowhere parallel to $\Theta$ on $M \setminus \pi^{-1}(\mathcal C)$, its pushforward $\pi_{*}R$ is a nowhere vanishing vector field on $\mathbb{S}^2 \setminus \mathcal C$. Moreover $\pi_{*}R$ is tangent to level sets of $K$, $K$ being perserved by the Reeb flow. Since $\gamma_i$ is transverse to the level sets of $K$, it is transverse to $\pi_*(R)$. Thus $R$ is transverse to $\Sigma = \pi^{-1}(\gamma_i)$. The level sets of $K$ in $\mathbb{S}^2 \setminus \mathcal C$ are unions of circles, hence the flow lines of $\pi_* R$ are circles as well. In particular, they pass infinitely many times through $\pi(x)$, which means that the Reeb orbit of $x$ in $M$ must intersect the $\mathbb{S}^1$-orbit of $x$, which is contained in $\Sigma_i$, infinitely many times. This proves the result. 
\end{proof}

For $x$ in $\Sigma$, we will write $\tau(x)$ for the first positive time at which the Reeb orbit starting at $x$ hits $\Sigma$ again. This allows us to define a diffeomorphism $\Phi$ of $\Sigma$ as

\[
	\begin{array}{ccccc}
		\Phi & \colon & \Sigma & \rightarrow & \Sigma \\
		& & x & \mapsto & \phi^{\tau(x)}(x)
	\end{array}
\]
as well as $\Phi_i = \Phi_{|\Sigma_i}$. Note that periodic points of $\Phi$ correspond to closed Reeb orbits of $\alpha$. 

\subsection{A decomposition for invariant contact forms} % (fold)
\label{sub:a_normal_form_for_invariant_contact_forms}

In this section, we give a decomposition for invariant contact forms on $M$. We start with a preliminary lemma.

\begin{lemma}
    \label{lemma:lift_form_quotient}
    Let $U$ be an $\mathbb{S}^1$-invariant open set of $M$, and $\kappa$ a one-form on $U$. Then there is a one-form $\beta$ on $\pi(U) \subset \mathbb{S}^2$ such that $\kappa = \pi^* \beta$ if and only if $\mathcal L_{\Theta} \kappa = 0$ and $i_{\Theta} \kappa = 0$.
\end{lemma}

\begin{proof}
	The direct implication is easy. Conversely, let $\kappa$ satisfying $\mathcal L_{\Theta} \kappa = 0$ and $i_{\Theta}\kappa = 0$. We define $\beta$ by $\beta_x(X) = \kappa_y(Y)$, where $\pi(y)=x$ and $\mathrm{d} \pi_y(Y) = X$. Since $i_{\Theta} \kappa = 0$, $\beta_y(Y)$ does not depend on the choice of $Y$. $\beta_y(Y)$ also does not depend on the choice of $y$ as $\mathcal L_{\Theta}\kappa = 0$.
\end{proof}

Since $\mathbb{S}^2 \setminus \pi(Z)$ is the union of two disks, and because the only $\mathbb{S}^1$-principal bundle over a disk is the trivial one, each of the two components of $M \setminus Z$ is a solid torus, the fundamental group of which being generated by the $\mathbb S^1$-orbits of the action. Let $\sigma$ be an $\mathbb S^1$-invariant closed one-form on $M \setminus Z$ such that $\int_{\mathbb{S}^1\cdot x} \sigma = 1$ for all $x$ in $M \setminus Z$, the orientation of $\mathbb S^1 \cdot x$ being given by the action. We will occasionally write $\sigma ^+$ (resp. $\sigma^-$) for the restriction of $\sigma$ to $(K>0)$ (resp. to ($K < 0$)).

Since the one-form $\alpha - K \sigma$ defined on $M \setminus Z$ is $\mathbb{S}^1$-invariant, and satisfies $i_\Theta(\alpha - K \sigma) = K - K = 0$, according to Lemma \ref{lemma:lift_form_quotient} we have on $M \setminus Z$,
$$
    \alpha = K \sigma + \pi^*\beta,
$$
for some one-form $\beta$ on $\mathbb{S}^2 \setminus \pi(Z)$.

We end this section with a useful geometrical interpretation of the Euler number $e$ of the $\mathbb{S}^1$-bundle $M \rightarrow \mathbb{S}^2$ (see for example \cite{jankins1983lectures} for a detailed discussion in the more general context of Seifert manifolds).

\begin{remark}
    \label{remark:euler_number}
    Let $T_1, T_2$ be two oriented solid tori, each of them being the total space of an $\mathbb{S}^1$-principal bundle $T_i \overset{\pi}{\rightarrow} D_i$ over a disk $D_i$. We write $h_i$ for an $\mathbb{S}^1$-fiber of $T_i \overset{\pi}{\rightarrow} D_i$, and $\mu_i$ for a circle in $T_i$ projecting bijectively to $\partial D_i$ (both being oriented respectively by the $\mathbb{S}^1$-action and by the natural boundary orientation of $\partial D_i$). If a bundle $M \rightarrow \mathbb{S}^2$ is obtained by gluing $T_1$ to $T_2$ via an attaching map $\phi \colon \partial T_1 \rightarrow \partial T_2$ whose action on fundamental groups writes $\phi_* = \begin{pmatrix} 1 & -e\\0 & -1\end{pmatrix}$ in the basis $(h_1, \mu_1)$ and $(h_2, \mu_2)$, then the integer $e$ coincides with the Euler number of $M$ defined earlier.
\end{remark}

We will write $c_x$ for the portion of Reeb orbit $\phi^{[0,\tau(x)]}$ starting at $x$ and ending at $\Phi(x)$. If $k \in ]K^-_i, K^+_i[$, there is a unique $x_{i,k}$ in $\tilde \gamma_i$ such that $K(x_{i,k}) = k$. We will occasionally write $c_{i,k}$ for $c_{x_{i,k}}$.

\begin{prop}
	\label{prop:equality_euler_number}
	The following identity holds:
	\[
		\int_{c_{i_0,0}} \sigma^+ - \sigma^- = e.
	\]
\end{prop}

\begin{proof}
    We write $D_1$ (resp. $D_2$) for the 2-disk $(K > 0)$ (resp. $(K<0)$), and $\text{Crit}(f)$ for the set of critical values of a real-valued function $f$. The bundle $M$ is obtained by gluing the two bundles $T_1 = \pi^{-1}(D_1)$ and $T_2 = \pi^{-1}(D_2)$ via the attaching map described above. To understand the attaching map, we must understand the orientations of $\partial T_1$ and $\partial T_2$. Let $Y$ be an $\mathbb{S}^1$-invariant gradient-like vector field for $\tilde K$ on $M$, i.e. a vector field vanishing exactly at critical points of $\tilde K$, and such that $\mathrm{d}\tilde K(Y) > 0$ otherwise. Since $\alpha \wedge \mathrm{d}\alpha (R,Y,\Theta) = \mathrm{d}\alpha(Y,\Theta) = \mathrm{d}K(Y) > 0$ on $M \setminus \text{Crit}(\tilde K)$, $(R,Y,\Theta)$ is a positive basis on $M \setminus \text{Crit}(\tilde K)$. Thus, $(\pi_* R, \pi_* Y)$ is a positive basis on $\mathbb{S}^2 \setminus \text{Crit}(K)$. We get that the induced boundary orientation of $\partial D_1$ (resp. of $\partial D_2$) coincides with (resp. is the opposite of) the orientation induced by the Reeb flow. Following Remark \ref{remark:euler_number}, the attaching map $\phi$ acts on the level of fundamental groups via the map $\phi_* = \begin{pmatrix} 1 & -e\\0 & -1\end{pmatrix}$ in the basis $(h_1, \mu_1)$ and $(h_2, \mu_2)$. Hence, the induced map $\phi^* \colon H^1(\partial T_2, \mathbb{Z}) \rightarrow H^1(\partial T_1, \mathbb{Z})$ is $\begin{pmatrix} 1 & 0\\-e & -1\end{pmatrix}$ in the dual basis $(h_1^*, \mu_1^*)$ and $(h_2^*, \mu_2^*)$. Since $\sigma^+ = h_1^*$ on $T_1$ and $\sigma^- = h_2^*$ on $T_2$, and since $\pi(c_{i_0,0}) = \mu_1$ in $\pi_1(\partial T_1)$, we get
    \[
        \int_{c_{i_0,0}} \sigma^+ - \sigma^- = \int_{c_{i_0,0}} h_1 - \phi^*{h_2} = \int_{\mu_1} e \mu_1^* = e.
    \]
\end{proof}

% subsection a_normal_form_for_invariant_contact_forms (end)

\subsection{A family of potentials} % (fold)
\label{sub:a_family_of_potentials}
For each $i \in I$ we define the map
\[
	\begin{array}{ccccc}
		J_i & \colon & ]K_i^-, K_i^+[ & \rightarrow & \mathbb{R} \\
		& & k & \mapsto & \int_{\pi(c_{i,k})} \beta.
	\end{array}.
\]

Since the curves $\pi(c_{i,k})$ are closed curves that uniformly converge to $\epsilon_i^\pm$ when $k$ goes to $K_i^\pm$, $J_i$ can be extended continuously to $[K^-_i, K^+_i]$. Moreover, since $\lim \limits_{k \rightarrow 0, k > 0} J_{i_0}(k)= \lim \limits_{k \rightarrow 0, k > 0} \int_{c_{i,k}} \alpha - k \sigma = \int_{c_{i_0,0}} \alpha$, and similarly $\lim \limits_{k \rightarrow 0, k < 0} J_{i_0}(k) = \int_{c_{i_0,0}} \alpha$, $J_{i_0}$ is continuous at $0$.

\begin{lemma}
    \label{lemma:derivative_J}
    The map $J_i$ is differentiable on $]K_i^-, K_i^+[ \setminus \{0\}$, with $J_i'(k) = -\int_{c_{i,k}} \sigma$. In addition, $$J'_{i_0}(0^+)-J'_{i_0}(0^-) = -e.$$
\end{lemma}

\begin{proof}
    Let $k^- < k^+$ be in $]K_i^-, K_i^+[$, such that $k^-$ and $k^+$ have same sign. Define $S$ to be $\bigcup \limits_{k^- < k < k^+}c_{i,k}$. This surface projects injectively on $\mathbb{S}^2$, with image $\Gamma_i \cap K^{-1}(]k^-,k^+[)$, whose oriented boundary is $\pi(k^-) - \pi(k^+)$ (both oriented by the flow of $\pi_* R$). By Stokes theorem, one then has
    \begin{align*}
        J_i(k^+)-J_i(k^-) &= -\int_{\pi(S)} \mathrm{d}\beta \\
        &= -\int_{S} \pi^*\mathrm{d}\beta \\
        &= -\int_{S} \mathrm{d}\alpha - \mathrm{d}K \wedge \sigma \\
        &= -\int_{S} \sigma \wedge \mathrm{d}K \\
        &= \int_{k^-}^{k^+}\left(-\int_{c_{i,k}} \sigma \right) \mathrm{d}k.
    \end{align*}
    The penultimate equality follows from the fact that the restriction of $\mathrm{d} \alpha$ to $TS$ vanishes as $R$ is tangent to $S$, and the last equality is obtained after the change of variable $x \in \tilde \gamma_i \cap S \mapsto K(x) \in ]k^-, k^+[$. Dividing by $k^+ - k^-$ and letting $k^+$ go to $k^-$ yields the result. Finally, $J'_{i_0}(0^+)-J'_{i_0}(0^-) = -e$ follows from Proposition \ref{prop:equality_euler_number}.
\end{proof}

An immediate corollary of Lemma \ref{lemma:derivative_J} is the identity below relating $\tau$ and $J_i$.

\begin{coro}
    \label{coro:identity_tau}
    For all $k \in ]K_i^-, K_i^+[$, $\tau(x_{i,k}) = J_i(k) - k J_i'(k)$.
\end{coro}

\begin{proof}
    Integrating the identity $\alpha = K \sigma + \pi ^* \beta$ along $c_{i,k}$ yields
    $$
        \tau(x_{i,k}) = k \int_{c_{i,k}}\sigma + J_i(k).
    $$
    By Lemma \ref{lemma:derivative_J}, $\int_{c_{i,k}}\sigma = - J_i'(k)$ and the result follows.
\end{proof}

A second corollary of Lemma \ref{lemma:derivative_J} is that the derivative of $J_i$ is closely related to periodic points of the first-return map of the surface of section.

\begin{coro}
    \label{coro:periodic_point}
    The Reeb orbit starting at $x_{i,k} \in \tilde \gamma_i$ is closed if and only if $J_i'(k)$ is rational. Moreover, if $J_i'(k) = \frac{p}{q}$ with $\gcd(p,q)=1$, then $x_{i,k}$ is a $q$-periodic point of $\Phi_i$, and the corresponding closed Reeb orbit has minimal period $q(J_i(k) - k J_i'(k))$.
\end{coro}

\begin{proof}
    Since $\tilde K$ is preserved by the Reeb flow, the first-return map $\Phi_i$ of $\Sigma_i$ induces a circle diffeomorphism of $\mathbb{S}^1 \cdot x_{i,k}$ for each $k \in ]K_i^-, K_i^+[$. By $\mathbb{S}^1$-invariance of the Reeb flow, each of those circle diffeomorphism is a translation. Finally, for $k \neq 0$, homotoping $c_{i,k}$ relative to its boundary in $M \setminus Z$ to a path contained in $\mathbb{S}^1 \cdot x_{i,k}$, together with the closedness of $\sigma$, gives us that the circle diffeomorphism of $\mathbb{S}^1 \cdot x_{i,k}$ is a translation of shift $\int_{c_{i,k}} \sigma$. As a consequence, the Reeb orbit starting at $x_{i,k}$ is closed if and only if $\int_{c_{i,k}} \sigma = -J_i(k)$ is a rational number.

    When $J_i'(k) = \frac{p}{q}$ with $\gcd(p,q)=1$, the translation of shift $\int_{c_{i,k}} \sigma = - \frac{p}{q}$ is $q$-periodic. Furthermore, the closed Reeb orbit starting at $x_{i,k}$ has period $q \tau(x_{i,k}) = q(J_i(k) - k J_i'(k))$ according to Corollary \ref{coro:identity_tau}.
\end{proof}

We now give some properties of the family of potentials $(J_i)_{i \in I}$. First of all, Lemma \ref{lemma:J_volume} tells us how the integral of the $J_i$'s are related to the volume.

\begin{lemma}
    \label{lemma:J_volume}
    The identity
    $$
        \vol(\alpha) = 2 \sum \limits_{i \in I} \int_{]K_i^-, K_i^+[}J_i \mathrm{d}k
    $$
    holds.
\end{lemma}

\begin{proof}
    Since $\sigma$ is closed, the identity $\alpha = K \sigma + \pi ^* \beta$ implies that
    $$
        \alpha \wedge \mathrm{d}\alpha = \sigma \wedge \pi^* (K \mathrm{d}\beta + \beta \wedge \mathrm{d}K).
    $$
    Let $T$ be a connected component of $M \setminus Z$. The volume of $T$ can be written as
    \begin{align*}
        \vol(T)&= \int_T \sigma \wedge \pi^* (K \mathrm{d}\beta + \beta \wedge \mathrm{d}K) \\
        &=\int_{\pi(T)} K \mathrm{d}\beta + \beta \wedge \mathrm{d}K \\
        &=\int_{\pi(T)} \mathrm{d}(K \beta) + 2 \beta \wedge \mathrm{d}K \\
        &=\int_{\pi(Z)}K \beta + 2\int_{\pi(T)} \beta \wedge \mathrm{d}K \\
        &=2\int_{\pi(T)} \beta \wedge \mathrm{d}K \\
        &=2\int_{\pi(T) \setminus \mathcal C} \beta \wedge \mathrm{d}K \\ 
        &=2\sum \limits_{i \in I} \int_{\Gamma_i \cap \pi(T)} \beta \wedge \mathrm{d}K \\
        &=2\sum \limits_{i \in I} \int_{]K_i^-, K_i^+[ \cap K(\pi(T))} \left(\int_{\pi(c_{i,k})}\beta\right) \mathrm{d}k.
    \end{align*}
    The last equality follows from the change of variable $x \in \gamma_i \cap \pi(T) \mapsto K(x) \in ]K_i^-, K_i^+[ \cap K(\pi(T))$. The claimed identity is then obtained by summing over both components $T$ of $M \setminus Z$.
\end{proof}

The key properties of the family $(J_i)_{i \in I}$ all follow from a differential inequality involving $J_i$.

\begin{prop}
    \label{prop:ineq_diff_J}
    For all $i$ in $I$, the map $J_i \colon ]K_i^-, K_i^+[ \setminus\{0\} \rightarrow \mathbb{R}$ satisfies 
    $$
        \forall k \in ]K_i^-, K_i^+[ \setminus\{0\}, J_i(k) - k  J_i'(k) > 0.
    $$
\end{prop}

\begin{proof}
    This follows immediately from Lemma \ref{lemma:derivative_J} and Corollary \ref{coro:identity_tau} in addition to the positivity of $\tau$.
\end{proof}

We now state some properties of this differential inequality.

\begin{lemma}
    \label{lemma:proprietes_J}
	Let $f \colon [a,b] \rightarrow \mathbb{R}$ be a continuous function, smooth on $]a,b[$, with $0 \notin ]a,b[$, such that for all $x \in ]a,b[$, 
    \begin{equation}
        f(x) - x f'(x) > 0.
        \label{eq:ineq_diff}
    \end{equation}
	If $f(b) \geq 0$ (resp. $f(a) \geq 0$) and $b > 0$ (resp. $a < 0$), then $f(x) > 0$ for all $x$ in $]a,b[$.
\end{lemma}

\begin{proof}
    Without loss of generality, we can assume that $f(b) \geq 0$ and $0 < a < b$. If there exists $x_0 \in ]a,b[$ such that $f(x_0) \leq 0$, it follows immediately from inequality (\ref{eq:ineq_diff}) at $x_0$, plus the fact that $x_0 > 0$, that $f'(x_0) < 0$. Furthermore, the inequality (\ref{eq:ineq_diff}) implies that $f(x) > 0$ whenever $f'(x) = 0$. Thus $f'$ is negative on $[x_0,b[$. As a consequence, $f(x) < 0$ on $[x_0, b]$ and in particular $f(b) < 0$, which is a contradiction.
\end{proof}

\begin{lemma}
    \label{lemma:int_ligne_niveau_positive}
    The function
    \[
        \begin{array}{ccccc}
            f & \colon & [\min(K), \max(K)] \setminus \text{Crit}(K) & \rightarrow & \mathbb{R} \\
            & & k & \mapsto & \int_{K^{-1}(k)}\beta
        \end{array}
    \]
    is positive.
\end{lemma}

\begin{proof}
    Let $T$ be a component of $M \setminus Z$. Assume without loss of generality that $\tilde K$ is positive on $T$.
    \begin{itemize}
        \item
            \textbf{When $K$ has isolated critical points: } Let $I_T = \{i \in I, \tilde \gamma_i \cap T \neq \emptyset\}$, and $i \in I_T$. We claim that when $K$ has isolated critical points, $f$ extends continuously to $[\min(K), \max(K)]$. Indeed, in that case if $k_0 \in \text{Crit}(K)$, the curve $K^{-1}(k)$ converges to $K^{-1}(k_0)$ for the Hausdorff topology when $k$ goes to $k_0$. As a consequence, $f(k)$ converges to $\int_{K^{-1}(k_0)} \beta$ when $k$ goes to $k_0$. Let us write $K_1 < \ldots < K_p = \max K$ for the ordered critical values of $K$ on $T$, and set $K_0 = 0$. $f$ is smooth on $]0, K_p] \setminus \{K_1, \ldots, K_p\}$, and continuous on $[0, K_p]$. Moreover, $f_{|]K_l, K_{l+1}[}$ is the sum of all the $J_i$ such that $]K_l, K_{l+1}[ \subset ]K_i^-, K_i^+[$. By Proposition \ref{prop:ineq_diff_J} applied to each of those $J_i$ and by linearity, $f$ satisfies the differential inequality
            \[
            	f(x) - x f'(x) > 0
            \]
            on each $]K_l, K_{l+1}[$. Since $K^{-1}(K_p)$ consists of isolated points, all of them maximizers for $K$, $f(K_p)$ is zero. Lemma \ref{lemma:proprietes_J} gives that $f$ is positive on $]K_{p-1}, K_p[$. Moreover, continuity of $f$ at every $K_i$ when $i < p$ allows us to apply Lemma \ref{lemma:proprietes_J} on each $]K_{p-l-1}, K_{p-l}[$ successively when $l$ goes from $1$ to $p-1$. This gives us the result.

        \item 
            \textbf{The general case: } Let $k \in [\min(K), \max(K)] \setminus \text{Crit}(K)$. Since $[\min(K), \max(K)] \setminus \text{Crit}(K)$ is open, let $[k-\eta, k + \eta]$ be a closed neighborhood of $k$ in $[\min(K), \max(K)] \setminus \text{Crit}(K)$, and let $\tilde g \colon M \rightarrow ]0, +\infty[$ a smooth function such that:
            \begin{itemize}
                \item $\tilde g$ is $\mathbb S^1$-invariant and hence induces a function $g \colon \mathbb S^2 \rightarrow ]0, + \infty[$;
                \item $\tilde g$ is constant equal to $1$ on $\tilde K^{-1}([k-\eta, k + \eta])$;
                \item $g K$ has isolated critical points.
            \end{itemize}
            Such a function $\tilde g$ exists since $[k-\eta, k + \eta]$ only contains regular values of $K$. Then $\tilde g \alpha$ is an $\mathbb S^1$-invariant contact form on $M$, with $\tilde K^{\tilde g \alpha} = \tilde g \tilde K$. Since $gK$ has isolated critical points, the previous argument implies that $f^{\tilde g \alpha}$ is positive on $[\min(gK), \max(gK)] \setminus \text{Crit}(gK)$. Since $\tilde g$ is constant equal to one on $\tilde K^{-1}([k-\eta, k + \eta])$, $K^{\tilde g \alpha} = K$ and $\beta^{ \tilde g \alpha} = \beta^\alpha$ on $K^{-1}([k-\eta, k + \eta])$, we have $f^{\tilde g \alpha} = f^\alpha$ on $[k-\eta, k + \eta]$. In particular, $f^\alpha(k) > 0$. The result then follows as $k$ is an arbitrary regular value of $K$.
    \end{itemize}
\end{proof}

\section{Proof of the main results}

\subsection{Proof of Theorem \ref{thm:ineg_sys_tight}}
In this section, we prove Theorem \ref{thm:ineg_sys_tight}. Since the case of $\mathbb S^1$-actions transverse to the contact structure has been taken care of in Proposition \ref{prop:ineg_sys_facile}, we focus now on the case where $e > 0$ and for which the set of Legendrian $\mathbb S^1$-orbits is non-empty and connected. We prove separately in Proposition \ref{prop:optimal_constant} the optimality of the constant $\frac{1}{2}$ in the systolic inequality when the Euler number is larger than $2$.

\begin{prop}
    \label{prop:optimal_constant}
    If $e > 2$, then the supremum of $\frac{\sys(\alpha)^2}{\vol(\alpha)}$ on $\Omega(e)$ is at least $\frac{1}{2}$.
\end{prop}

\begin{proof}
    We will construct contact forms $\alpha$ in $\Omega(e)$ such that $\frac{\sys(\alpha)^2}{\vol(\alpha)}$ is arbitrarily close to $\frac{1}{2}$. We fix $\eta \in ]0, \frac{1}{2+e}[$, $a = \frac{1}{2} - \frac{e \eta}{2} > 0$ and $f \colon ]-\frac{\eta}{2}, a[ \rightarrow \mathbb{R}$, a smooth function satisfying:
    \begin{itemize}
        \item for all $r \in ]-\frac{\eta}{2}, \frac{\eta}{2}[$, $f(r) = \frac{1}{2} - \frac{e}{2}r$;
        \item for all $r \in ]\eta, a[$, $f(r) = a - r$;
        \item $f$ is convex on $]\frac{\eta}{2}, \eta[$.
    \end{itemize}

    See Figure \ref{fig:interpolation} for an example of such a function. Taking $f$ to be a convex smoothening of the piecewise linear function defined by
    \[
        r \mapsto
        \begin{cases}
            \frac{1}{2} - \frac{e}{2}r &\text{ if } - \frac{\eta}{2} < r \leq \frac{e \eta}{e+2} \\ 
            a - r &\text{ if } \frac{e \eta}{e+2} \leq r < a
        \end{cases}
    \]
    gives a function as requested since $\frac{\eta}{2} < \frac{e \eta}{e+2} < \eta$ when $e > 2$.

    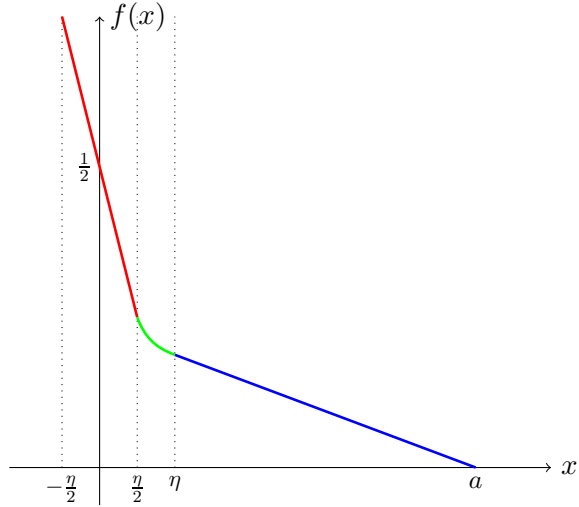
\begin{figure}[h]
        \centering
        \begin{tikzpicture}[scale=1]
            \draw[->] (-1.2,0) -- (6,0) node[right] {$x$};
            \draw[->] (0,-0.5) -- (0,6) node[right] {$f(x)$};
            \draw[dotted] (0.5,0) -- (0.5,6);
            \draw[dotted] (1,0) -- (1,6);
            \draw[dotted] (-0.5,0) -- (-0.5,6);
            \draw[]
                node[below, scale=0.8] at (-0.5,0) {$-\frac{\eta}{2}$}
                node[below, scale=0.8] at (0.5,0) {$\frac{\eta}{2}$}
                node[below, scale=0.8] at (1,0) {$\eta$}
                node[below, scale=0.9] at (5,0) {$a$}
                node[left, scale=0.8] at (0,4) {$\frac{1}{2}$};
            \draw[line width=1pt, color=red] (-0.5,6) -- (0.5,2);
            \draw[line width=1pt, color=blue] (1,1.5) -- (5,0);
            \draw[line width=1pt, color=green] (0.5,2) to[out=-70, in=160] (1,1.5);
        \end{tikzpicture}\\
        \label{fig:interpolation}
        \caption{\small Example of a smooth convex function $f$ defined on $]-\frac{\eta}{2}, a[$, satisfying $f(r)=\frac{1}{2} - \frac{e}{2}r$ on $]-\frac{\eta}{2}, \frac{\eta}{2}[$ and $f(r)=a-r$ on $]\eta, a[$.}
    \end{figure}

    We now define two contact forms $\alpha_+$ and $\alpha_-$ on the solid torus $T = \mathbb{S}^1 \times D^2$ as follows:
    \[
        \alpha_+ = (a-r^2)\mathrm{d}t + f(a-r^2) \mathrm{d}s
    \]
    and
    \[
        \alpha_- = -(a-r^2)\mathrm{d}t - f(a-r^2) \mathrm{d}s
    \]
    where $t$ is the coordinate along the $\mathbb{S}^1$ factor and $\left(r,s\right) \in \left[0, \sqrt{a+\frac{\eta}{2}}\right] \times \mathbb{S}^1$ are polar coordinates on the 2-disk of radius $\sqrt{a+ \frac{\eta}{2}}$. Remark that $\alpha_+$ and $\alpha_-$ are indeed contact forms since the contact condition $\alpha_\pm \wedge \mathrm{d} \alpha_\pm  > 0$ becomes the condition $f(a-r^2) - (a-r^2) f'(a-r^2) >0$, which is easily seen to be satisfied. We write $U$ for the neighborhood $\left]\sqrt{a-\frac{\eta}{2}}, \sqrt{a+\frac{\eta}{2}}\right] \times \partial T$ of $\partial T$. Let $g \colon U \rightarrow U$ be the diffeomorphism defined by $g(r,s,t) = (2a-r^2,-(s+et),t)$. Following Remark \ref{remark:euler_number}, gluing $T$ to itself along $U$ via $f$ produces the manifold $M$ that is the total space of the $\mathbb{S}^1$-principal bundle over $\mathbb{S}^2$ with Euler number $e$. Furthermore, since 
    \[
        (g^* \alpha_+)_{r,s,t} = -(a-r^2) \mathrm{d}t - (f(r^2-a) + e(r^2-a)) \mathrm{d}s
    \]
    and since $f(r^2-a) + e (r^2-a) = f(a-r^2)$ when $(r^2-a) \in ]- \frac{\eta}{2}, \frac{\eta}{2}[$, we obtain that
    \[
        (g^* \alpha_+)_{r,s,t} = - (a-r^2) \mathrm{d}t - f(a-r^2) \mathrm{d}s = (\alpha_-)_{r,s,t}
    \]
    holds true on $U$. As a consequence, gluing $\alpha_+$ and $\alpha_-$ via $g$ produces a contact form $\alpha \in \Omega(e)$. This invariant contact form satisfies the following properties :
    \begin{itemize}
        \item Its moment map $K^\alpha(r,s,t) = \iota_{\partial t} \alpha_{\pm}(r,s,t) = \pm (a-r^2)$ has image the interval $[-a,a]$, and has exactly two critical points, with critical values $\pm a$. They correspond to two closed Reeb orbits that coincide with orbits of the $\mathbb{S}^1$-action, which have both period $a$.
        \item Since $K^\alpha$ only has two critical points, the family of surfaces of sections constructed in Section \ref{sub:invariant_surface_of_section} consists in only one surface of section, corresponding to the interval $]-a,a[$ of regular values of $K^\alpha$.
        \item On $(K^\alpha>0)$, the decomposition of $\alpha$ described in Section \ref{sub:a_normal_form_for_invariant_contact_forms} is given by the identity $\alpha_+ = (a-r^2)\mathrm{d}t + f(a-r^2) \mathrm{d}s$. With the notations of Section \ref{sub:a_normal_form_for_invariant_contact_forms}, one has $\sigma = \mathrm{d}t$, $K(r,s,t)=a-r^2$ and $\beta_{r,s,t} = f(a-r^2)\mathrm{d}s$. Similarly, one has $\sigma = \mathrm{d}t$, $K(r,s,t)=-(a-r^2)$ and $\beta_{r,s,t} = -f(a-r^2)\mathrm{d}s$ on $(K^\alpha < 0)$.
        \item Since there is only one surface of section, we omit the index and write $J^\alpha \colon ]-a, a[ \rightarrow \mathbb{R}$ for the corresponding potential. It follows from the identity $\beta_{r,s,t} = f(a-r^2)\mathrm{d}s$ on $(K^\alpha>0)$ that $J^\alpha(k) = f(k)$ for $k \in [0,a[$. Similarly, $J^\alpha(k) = f(-k)$ for $k \in ]-a,0]$. 
    \end{itemize}
    Except for the two closed Reeb orbits that coincide with $\mathbb{S}^1$-orbits, the closed Reeb orbits of $\alpha$ are in correspondance with points where $(J^{\alpha})'$ is rational, with period given by Corollary \ref{coro:periodic_point}.
        \begin{itemize}
            \item If $|r| > \eta$, then $|(J^{\alpha})'(r)| = |f'(r)| = 1$, and the corresponding closed Reeb orbit has period $a$;
            \item if $|r| < \frac{\eta}{2}$, then $|(J^{\alpha})'(r)| = |f'(r)| = \frac{e}{2}$, and the corresponding closed Reeb orbit has period $\frac{1}{2}$ if $e$ is even, or $1$ if $e$ is odd;
            \item if $|r| \in ]\frac{\eta}{2},\eta[$ and $|(J^{\alpha})'(r)| = |f'(r)| = \frac{p}{q}$, then the corresponding closed Reeb orbit has period 
            \begin{align*}
                q(J^{\alpha}(r) - r (J^{\alpha})'(r)) &= q(f(|r|) - r f'(|r|))\\
                & \geq f(|r|) - r f'(|r|) \\ 
                & \geq a  
            \end{align*}
            where the last inequality follows from the convexity of $f$.
        \end{itemize}
        As a consequence, $\sys(\alpha)=a$, which converges to $\frac{1}{2}$ when $\eta$ converges to zero. Moreover, Lemma \ref{lemma:J_volume} implies that $\vol(\alpha)$ converge to $\frac{1}{2}$ when $\eta$ converges to zero. We conclude that $\lim \limits_{\eta \rightarrow 0}\frac{\sys(\alpha)^2}{\vol(\alpha)} = \frac{1}{2}$, which proves the result.
\end{proof}

\begin{remark}
    The contact form $\alpha$ constructed in the proof of Proposition \ref{prop:optimal_constant} is briefly described in \cite[Section 3]{MR4529024}, as a deformation of the round metric on $\mathbb{S}^2$ to a spindle orbifold.
\end{remark}

We can now prove Theorem \ref{thm:ineg_sys_tight}.

\begin{proof}{(Theorem \ref{thm:ineg_sys_tight})}
    Let $\alpha \in \Omega(e)$, where $e > 0$. We write $J$ for $J_{i_0}$, and $K_{\text{min}}$ for the minimal absolute value of a critical value of $K$. Since $]-K_{\text{min}}, K_{\text{min}}[$ is contained in $U_{i_0}$ and does not intersect $]K_i^-, K_i^+[$ when $i \neq i_0$, Lemma \ref{lemma:int_ligne_niveau_positive} tells us that $J$ is positive on $[-K_{\text{min}}, K_{\text{min}}]$. By Lemma \ref{lemma:J_volume} and \ref{lemma:int_ligne_niveau_positive}, we have 
	\begin{equation}
		\label{eq:volume_J_degre_1}
		\vol(\alpha) \geq 2 \int_{-K_{\text{min}}}^{K_{\text{min}}} J \mathrm{d}k.
	\end{equation}

    We first prove that $\sys(\alpha)^2 \leq \frac{1}{2} \vol(\alpha)$ when $e \geq 2$. In that case, we write $g_e$ for the function $k \mapsto J(k) + |k|$. By virtue of Lemma \ref{lemma:derivative_J}, $g_2$ is continuously differentiable on $]-K_{\text{min}}, K_{\text{min}}[$, and $g_e$ is continuously differentiable on $]-K_{\text{min}}, K_{\text{min}}[ \setminus \{0\}$. We claim that if $e \geq 2$,
    \begin{equation}
		\label{eq:lower_bound_g_e}
        g_e(k) \geq \max\left(|k|, \sys(\alpha)\right)
	\end{equation}
    Indeed, since $J$ is positive, one has $g_e(k) \geq |k|$. Furthermore, as a continuous function, $g_e$ reaches its minimum on $[-K_{\text{min}}, K_{\text{min}}]$ at a point $k_0$. Three cases can occur : 
    \begin{itemize}
        \item If $k_0 \in \{-K_{\text{min}}, K_{\text{min}}\}$, then $g_e(k_0) \geq K_{\text{min}} \geq \sys(\alpha)$;
        \item If $k_0 \in ]-K_{\text{min}}, K_{\text{min}}[ \setminus \{0\}$, then $k_0$ is a critical point of $g_e$. Thus we have $|J'(k_0)| = 1$. By Corollary \ref{coro:periodic_point}, $k_0$ corresponds to a fixed point of the first-return map of the surface of section, hence a periodic orbit of the Reeb flow, of period the corresponding critical value of $g_e$. In particular, $g_e(k_0) \geq \sys(\alpha)$.
        \item If $k_0 = 0$, then for $g_e$ to be minimal at 0, we must have that $J(k) + k \geq J(0)$ for $k$ positive, implying $J'(0^+) \geq -1$. Similarly, we must have $J(k) - k \geq J(0)$ for $k$ negative, which gives us $J'(0^-) \leq 1$. Since $J'(0^+)-J'(0^-) = -e$, this can only happen when $e = 2$, in which case we have $J'(0^+) = -1$ and $J'(0^-) = 1$. As in the previous case, Corollary \ref{coro:periodic_point} gives that the corresponding Reeb orbit is closed, with period $g_e(k_0)$, which is larger than $\sys(\alpha)$.
    \end{itemize}

    We get that $\min(g_e) \geq \sys(\alpha)$, which together with $g_e(k) \geq |k|$ gives the lower bound (\ref{eq:lower_bound_g_e}). Additionally, the following holds:
    \begin{equation}
        \begin{split}
    		\label{eq:integrale_J_e}
    		\int_{-K_{\text{min}}}^{K_{\text{min}}} J \mathrm{d}k &= \int_{-K_{\text{min}}}^{K_{\text{min}}} g_e(k) \mathrm{d}k - \int_{-K_{\text{min}}}^{K_{\text{min}}}|k| \mathrm{d}k \\
    		&= \int_{-K_{\text{min}}}^{K_{\text{min}}} g_e(k) \mathrm{d}k - K_{\text{min}}^2 \\
            &\geq \int_{-K_{\text{min}}}^{K_{\text{min}}} \max\left(|k|, \sys(\alpha)\right) \mathrm{d}k - K_{\text{min}}^2 \\
            &= \sys(\alpha)^2 + K_{\text{min}}^2 - K_{\text{min}}^2 \\
            &= \sys(\alpha)^2.
    	\end{split}
    \end{equation}
    The two lower bounds (\ref{eq:volume_J_degre_1}) and (\ref{eq:integrale_J_e}) together yield the systolic inequality for $e\geq 2$, namely that

    \begin{equation}
        \label{eq:ineg_sys_e_2}
        \sys(\alpha)^2 \leq \frac{1}{2} \vol(\alpha).
    \end{equation}

    When $e=1$, we define $g_1$ as the function $k \mapsto J(k) + \frac{|k|}{2}$. The lower bound
    \begin{equation}
		\label{eq:lower_bound_g_1}
        g_1(k) \geq \max\left(\frac{|k|}{2}, \frac{\sys(\alpha)}{2}\right)
	\end{equation}
    holds for similar reasons to (\ref{eq:lower_bound_g_e}). Furthermore, plugging the lower bound \ref{eq:lower_bound_g_1} in the equation (\ref{eq:volume_J_degre_1}) gives
    \begin{equation}
        \begin{split}
    		\label{eq:integrale_J_1}
            \vol(\alpha) &\geq 2 \int_{-K_{\text{min}}}^{K_{\text{min}}} J \mathrm{d}k \\
            &= 2\int_{-K_{\text{min}}}^{K_{\text{min}}} g_1(k) \mathrm{d}k - 2\int_{-K_{\text{min}}}^{K_{\text{min}}}|k| \mathrm{d}k \\
            &\geq \int_{-K_{\text{min}}}^{K_{\text{min}}} \max\left(|k|, \sys(\alpha)\right) \mathrm{d}k - K_{\text{min}}^2 \\
            &= \sys(\alpha)^2 + K_{\text{min}}^2 - K_{\text{min}}^2 \\
            &= \sys(\alpha)^2.
    	\end{split}
    \end{equation}
    This proves the systolic inequality when $e=1$, namely

    \begin{equation}
        \label{eq:ineg_sys_e_1}
        \sys(\alpha)^2 \leq \vol(\alpha).
    \end{equation}

    We now study the equality case. When $e \geq 2$ (resp. $e=1$), the inequality (\ref{eq:ineg_sys_e_2}) is an equality if and only equality holds in (\ref{eq:volume_J_degre_1}) and (\ref{eq:lower_bound_g_e}) (resp. in (\ref{eq:volume_J_degre_1}) and (\ref{eq:lower_bound_g_1})). Equality holds in (\ref{eq:volume_J_degre_1}) if and only if $K$ has exactly two critical values $K_{\text{min}}$ and $-K_{\text{min}}$. In particular, all Reeb orbits that do not intersect the surface of section $\Sigma$ are closed, with period $K_{\text{min}}$. Furthermore, equality holds in (\ref{eq:lower_bound_g_e}) if and only if $J(k) = \sys(\alpha) - |k| = K_{\text{min}} - |k|$, which only happens if $e=2$ by Lemma \ref{lemma:derivative_J}. In that case, $|J'(k)| = 1$ for all $k$, thus $x_{i_0,k}$ is a fixed point of the surface of section. The Reeb orbit starting at any of the $x_{i_0,k}$ is then closed, with minimal period $\tau(k) = J(k) - k J'(k) = K_{\text{min}}$.
    
    Similarly, when $e=1$, equality holds in (\ref{eq:lower_bound_g_1}) if and only if $J(k) = \frac{\sys(\alpha)}{2} - \frac{|k|}{2} = \frac{K_{\text{min}}}{2} - \frac{|k|}{2}$. In that case, $|J'(k)| = \frac{1}{2}$ for all $k$, thus $x_{i_0,k}$ is a 2-periodic point of the surface of section. The Reeb orbit starting at any of the $x_{i_0,k}$ is then closed, with minimal period $2\tau(k) = 2(J(k) - k J'(k)) = K_{\text{min}}$.

    In both cases, we get that every Reeb orbit of $\alpha$ is closed, with common minimal period $K_{\text{min}}$. In other words, $\alpha$ is Zoll. It follows from the previous argument that the equality case cannot be reached when $e>2$. Finally, Proposition \ref{prop:optimal_constant} implies that the systolic inequality is sharp also when $e > 2$.
\end{proof}

\subsection{Proof of Corollary \ref{coro:sys_contr}}
\label{sub:proof_coro_contractible}

We define the prime action spectrum of a contact form $\alpha$ on a manifold $M$, denoted by $\sigma_P(M,\alpha)$, as the set of prime periods of closed Reeb orbits of $\alpha$. We will make use of the following theorem of Cristofaro-Gardiner and Mazzucchelli (Theorem 1.5 in \cite{MR4152621}), stating that two Besse contact forms on a closed three-dimensional manifold are diffeomorphic if and only if they share the same prime action spectrum.

\begin{thm}[Theorem 1.5 in \cite{MR4152621}]
    \label{thm:Besse_CGM}
    Let $M$ a closed connected $3$-manifold, and $\alpha_1, \alpha_2$ two Besse contact forms on $M$. Then $\sigma_P(M,\alpha_1) = \sigma_P(M,\alpha_2)$ if and only if there exists a diffeomorphism $\Psi \colon M \rightarrow M$ such that $\Psi^* \alpha_2 = \alpha_1$.
\end{thm}

We now give the proof of Corollary \ref{coro:sys_contr}.

\begin{proof}{(Corollary \ref{coro:sys_contr})}
    When $e \neq 0$, $M$ has fundamental group $\mathbb{Z}/|e| \mathbb{Z}$, which is torsion. As a consequence, every closed Reeb orbit of $\alpha \in \Omega(e)$ has an iterate which is contractible. The classification of $\mathbb S^1$-principal bundles over $\mathbb S^2$ yields that the circle action on $M$ lifts to a free circle action on $\mathbb S^3$ via any degree $|e|$ covering map $\mathbb S^3 \overset{p}{\rightarrow} M$. Moreover, the contact form $p^* \alpha$ on $\mathbb{S}^3$ is invariant under the lifted circle action. We observe that contractible closed Reeb orbits of $\alpha$ are precisely the closed orbits which lift to closed Reeb orbits of $p^* \alpha$, meaning that a Reeb orbit $\tilde c$ of $p^* \alpha$ with basepoint $x \in \mathbb{S}^3$ is closed if and only if the Reeb orbit $c$ of $\alpha$ starting at $p(x)$ is closed, and in that case the period of $\tilde c$ is the shortest period of a contractible iterate of $c$. It follows that the systole of $p^* \alpha$ agrees with the contractible systole of $\alpha$. Combining this discussion with the fact that the covering map $\mathbb{S}^3 \overset{p}{\rightarrow} M$ has degree $|e|$, we get that Theorem \ref{thm:ineg_sys_tight} for $|e| = 1$ implies the inequality
    \[
        \sys_{\text{contr}}(\alpha)^2 \leq |e| \vol(\alpha).
    \]
    Equality holds if and only if the lifted form $p^* \alpha$ is Zoll. Since the contractible systolic ratio $\frac{\sys_{\text{contr}}(\alpha)^2}{\vol(\alpha)}$ does not depend on $p$, whether or not $p^* \alpha$ is Zoll does not depend on the specific choice of covering map $p$. If $p^*\alpha$ is Zoll, then $\alpha$ is Besse since the Reeb orbits of $p^* \alpha$, which are closed, project on $M$ to the Reeb orbits of $\alpha$. However not all Besse forms on $M$ can be lifted to Zoll forms on $\mathbb S^3$. We spend the rest of the proof on proving a characterization of such Besse forms on $M$.

    \begin{itemize}
        \item When $e < 0$: if $p^* \alpha$ is Zoll, it follows from Proposition \ref{prop:ineg_sys_facile} that the lifted circle action on $\mathbb S^3$ coincides with the Reeb flow of $p^* \alpha$. Therefore the free circle action on $M$ coincides with the Reeb flow of $\alpha$, which is hence Zoll. Conversely, the Reeb flow of a Zoll contact form in $\Omega(e)$ coincides with the $\mathbb{S}^1$-action, which lifts to a free circle action on $\mathbb{S}^3$ via any degree $-e$ covering map. Hence Zoll contact forms in $\Omega(e)$ lift to Zoll contact forms on $\mathbb{S}^3$.

        \item When $e > 0$: if $p^* \alpha$ is Zoll, the proof of Theorem \ref{thm:ineg_sys_tight} implies that $K^{p^* \alpha}$ has exactly two critical points, meaning that exactly two Reeb orbits of $p^* \alpha$ coincide with orbits of the circle action. If $T>0$ is the minimal period of the Reeb flow of $p^* \alpha$, those two orbits correspond to two closed Reeb orbits of $\alpha$ of minimal period $\frac{T}{e}$. All the other orbits of $p^* \alpha$ project to closed Reeb orbits of $\alpha$ of period $T$ when $e$ is odd, $\frac{T}{2}$ when $e$ is even.  Thus if $e \in \{1,2\}$, $\alpha$ is Zoll as well, but otherwise $\alpha$ is only Besse and has two singular orbits, which have the same minimal period.

        Conversely, Zoll contact forms in $\Omega(e)$ lift to Zoll contact forms on $\mathbb{S}^3$ when $e \in \{1,2\}$. Assume now that $e>2$. We claim that if $\alpha$ is a Besse contact form on $M$ with exactly two singular orbits, which have the same minimal period, then $\alpha$ lifts to a Zoll contact form on $\mathbb{S}^3$.

        Indeed, the Reeb flow of any Besse contact form induces a Seifert fibration of $M$, namely a circle action on $M$ with finite isotropy groups. For any coprime integers $k_1, k_2 \in \mathbb{Z} \setminus \{0\}$, we define a $\mathbb{S}^1$-action on $M$, whose orbits define a Seifert fibration of $M$, as follows : the $\mathbb{S}^1$-action on $\mathbb{S}^3$ defined by
        \[
            t \cdot (z_1,z_2) \in \mathbb{S}^3 \mapsto (e^{i k_1 t} z_1, e^{i k_2 t} z_2)
        \]
        commutes with the $\mathbb{Z} / e \mathbb{Z}$ action generated by $(z_1, z_2) \mapsto (e^{\frac{2 \pi i}{e}} z_1, e^{-\frac{2 \pi i}{e}} z_2)$. As such it induces a circle action on the quotient $\mathbb{S}^3 \overset{p_0}{\rightarrow} M$ by the latter action. When $(k_1,k_2) = (1, -1)$ or $(-1,1)$, this defines a free action on $M$. Moreover, when $(k_1,k_2) = (1, 1)$ or $(-1,-1)$, the circle action on $\mathbb{S}^3$ coincides (up to orientation) with the Reeb flow of the standard contact form $\alpha_0 = \frac{1}{4i}\sum\limits_{i=1}^{2}\left(\bar{z_i} \mathrm{d}z_i - z_i \mathrm{d}\bar{z_i} \right)$, which is Zoll. $\alpha_0$ induces a Besse form $\alpha_1$ on $M$, with prime action spectrum $\left(\frac{T}{e}, T\right)$ if $e$ is odd, $\left(\frac{T}{e}, \frac{T}{2}\right)$ if $e$ is even. The circle action on $M$ for $(k_1,k_2) = (1, 1)$ or $(-1,-1)$ coincides with the Reeb flow of $\alpha_1$. Finally, the regular orbits of the circle action on $M$ for $(k_1,k_2) = (1, -1)$ or $(-1,1)$ have period $1$ if $e$ is odd, $\frac{1}{2}$ if $e$ is even. It has two singular orbits, which have period $\frac{1}{e}$.

        As a consequence of the classification of Seifert fibrations on lens spaces (see \cite[Theorem 5.1]{Geiges2018}), there exist coprime integers $k_1, k_2 \in \mathbb{Z} \setminus \{0\}$ and a constant $T > 0$ such that the Reeb flow of $\frac{\alpha}{T}$ is conjugated via a diffeomorphism of $M$ to the $\mathbb{S}^1$-action on $M$ associated to $(k_1,k_2)$. For the two singular orbits to have the same minimal period, $|k_1|$ and $|k_2|$ must be equal. Since $k_1$ and $k_2$ are coprime, they must belong to $\{-1,1\}$. If $k_1 = -k_2$, the induced action on $M$ is free, hence $\alpha$ is Zoll. The proof of Theorem \ref{thm:ineg_sys_tight} implies that this cannot happen when $e>2$. Consequently $k_1 = k_2$, and in that case the prime action spectrum of $\alpha$ is $\left(\frac{T}{e}, T\right)$ if $e$ is odd, $\left(\frac{T}{e}, \frac{T}{2}\right)$ if $e$ is even. It follows that $\frac{\alpha}{T}$ has the same prime action spectrum as $\alpha_1$. Applying Theorem \ref{thm:Besse_CGM}, we get that there is a diffeomorphism $\psi$ of $M$ such that $\psi^*\alpha=\alpha_1$. Finally, this gives that $(\psi \circ p_0)^* \alpha = p_0^* \alpha_1$ is a Zoll contact form on $\mathbb S^3$. As a consequence, when $e> 2$ a Besse contact form on $M$ lifts to a Zoll contact form on $\mathbb S^3$ if and only if it has two singular orbits, which have the same prime period. This conclude the proof.
    \end{itemize}
\end{proof}

\bibliographystyle{plain}
\bibliography{biblio}

\end{document}